\documentclass[12pt,oneside]{article}
\usepackage{amsmath,amsthm,amsfonts,amssymb,MnSymbol,mathrsfs,latexsym}
\usepackage{mathtools}
\usepackage[english]{babel}

\newcommand{\const}{\rm const}

\DeclareMathOperator*{\esssup}{ess\,sup}

\theoremstyle{plain}

\title{\large \textbf{Moment and exponential estimation for the distribution\\ of the norms for random matrices martingales}}

\footnotesize\date{}

\author{\normalsize Maria Rosaria Formica ${}^{1}$,   \normalsize Eugeny Ostrovsky
${}^2$ and \normalsize Leonid Sirota ${}^3$}

\begin{document}

  \maketitle

\begin{center}
{\footnotesize ${}^{1}$ Universit\`{a} degli Studi di Napoli \lq\lq Parthenope\rq\rq, via Generale Parisi 13,\\
Palazzo Pacanowsky, 80132,
Napoli, Italy.} \\

\vspace{1mm}

{\footnotesize e-mail: mara.formica@uniparthenope.it} \\

\vspace{2mm}

{\footnotesize ${}^{2,\, 3}$  Bar-Ilan University, Department of Mathematics and Statistics, \\
52900, Ramat Gan, Israel.} \\

\vspace{1mm}

{\footnotesize e-mail: eugostrovsky@list.ru}\\

\vspace{1mm}

{\footnotesize e-mail: sirota3@bezeqint.net} \\

\end{center}

\vspace{3mm}

\begin{abstract}
 We derive sharp non - asymptotical Lebesgue - Riesz as well as Grand Lebesgue Space norm estimations for different  norms of matrix martingales through these norms for the correspondent martingale differences and through the entropic dimension of the extremal points of the unit ball for a basic space.\par
 \ These estimates allow us to deduce in particular the exponential decreasing tail of distribution for these norms of matrix martingales. \par
 \ We bring also some examples in order to show the exactness of the obtained  estimations.\par

\end{abstract}

\vspace{3mm}


 \vspace{5mm}

 \section{Introduction.}

 \vspace{5mm}

 \hspace{3mm} Let $ (\Omega = \{\omega\}, \cal{B}, {\bf P}) $ be a probability space (sufficiently rich)  with expectation $ \ {\bf E} \ $. The classical Lebesgue - Riesz norms for numerical (or complex) valued random variables (r.v.) $\xi$ are defined, as usual, by

$$
|\xi|_p \stackrel{def}{=} \ \left(\ {\bf E} |\xi|^p \ \right)^{1/p}, \ p \in [1,\infty);
$$

$$
|\xi|_{\infty} \stackrel{def}{=} \esssup_{\omega \in \Omega} |\xi(\omega)|.
$$

 \ Let  $ \ \{F_k\}, \ k = 0, 1,2,\ldots \ $ be certain {\it filtration}, i.e. be an increasing sequence of
sigma - subfields of the source sigma field $ \ \cal{B}:  $

$$
  F_0 = \{ \ \emptyset, \Omega \ \}, \hspace{3mm}  F_k \ \subset F_{k + 1} \subset \cal{B}.
$$

\hspace{3mm} Let also $ \ V(k), k = 0,1,2,\ldots,n;  \ n \le \infty \ $ be a   {\it matrix valued }
martingale of a size   $ \ d \otimes d \ $ with entries correspondingly

$$
V(k) = \{  \ V_{i,j}(k) \ \}, \hspace{3mm} i,j = 1,2,\ldots,d;
$$

\begin{equation} \label{matr mart}
V_0 = 0, \hspace{3mm}  {\bf E} V(k)/F_{k-1} = V(k-1), \ k = 1,2,3,\ldots.
\end{equation}

 \ Introduce the correspondent matrix valued {\it martingale differences,}  briefly  (mart dif) random variables

$$
\xi_{i,j}(k) := V_{i,j}(k) - V_{i,j}(k-1); \ k > 1; \  \xi_{i,j}(0) = 0;
$$
so that

\begin{equation} \label{sum represent}
V_{i,j}(n) = \sum_{l=1}^n \xi_{i,j}(l), \hspace{3mm} n = 1,2,\ldots \ .
\end{equation}

\vspace{3mm}

\hspace{3mm} Define  also the following important for us variables, Lebesgue - Riesz norms for the martingale differences
and for the values $ \ p \ge 1 \ $

 \vspace{3mm}

 \begin{equation} \label{mom mart dif}
 \mu_{i,j}(k)[p] =  \mu_{i,j}^{(n)}(k)[p] \stackrel{def}{=} |\xi_{i,j}(k)|_p, \hspace{3mm} \mu(k)[p] := \mu_{1,1}(k)[p];
 \end{equation}

\vspace{3mm}

\begin{equation} \label{mom mart dif}
\nu_{i,j}[p] = \nu_{i,j}^{(n)}[p] \stackrel{def}{=} \max_{k = 1,2, \ldots,n} |\xi_{i,j}|_p, \hspace{4mm}  \ \nu[p]:= \nu_{1,1}[p].
\end{equation}

\vspace{4mm}

\hspace{3mm} {\bf Our target in this preprint is  a tail estimate for the arbitrary norm of random
martingale matrix through the correspondent estimate for its martingale differences and the metric
 entropy of the extremal points for the unit ball of the correspondent finite - dimensional space. } \par

\vspace{4mm}

 \ Offered here estimations generalized ones in the works \cite{Mackey},  \cite{Tian}, \cite{Tropp1} - \cite{Tropp4} etc.
The so - called  "continuous time" matrix martingales  is studied in  \cite{Bacry}. \par

\vspace{4mm}

\section{One dimensional case. \ Khintchine statement.}

\vspace{3mm}

 \hspace{3mm} In this rubric $ \ d = 1, \  \{\xi(i)\}, \ i = 1,2,\ldots,n \ $ be a sequence of centered  mart - dif
relative the source  filtration; $ \ V(n) = \sum_{i=1}^n \xi(i) \ $ be a martingale  with correspondent mart - dif sequence $ \ \{\xi(i) \}. \ $ \par
 \hspace{3mm} Let also \ $ \ b = \vec{b} \ $ be a deterministic numerical sequence such that

$$
b = \vec{b} = \{b(1), b(2), \ldots, b(n)\}, \ \vec{b} \in B \ \Leftrightarrow \sum_{i=1}^n b^2(i) = 1,
$$

$$
Q(n, b; \xi(\cdot)) \stackrel{def}{=} \sum_{i=1}^n b(i) \ \xi(i),
$$

$$
  U = U(p; n) := \sup_{b \in B} |Q(n,b; \xi(\cdot))|_p,
$$
 the so - called {\it Khintchine variable.} \par
 \  We intent in this section to evaluate the variable $ \ U = U(p,n), \ p,n \ge 1 \ $ thought the moment
 values  $ \   \{ \ \mu(k)[p] \}, \ $ following via the values $ \ \nu(p); \ $  the so - called Khintchine statement of problem. \par
 \ There are a huge of works devoting this problem in the classical statement, see for instance
\cite{Burkholder1}, \cite{Burkholder2}, \cite{Fan}, \cite{Freedman}, \cite{Hitchenko}, \cite{Lesign Volny},
\cite{Liu}, \cite{Pena}, \cite{Peshkir  Shirjaev}, \cite{Ratchkauskas}, \cite{Volny2}. \par

\vspace{3mm}

 \  Let us take into account the following important variable belonging  to Osekowski A, see \ \cite{Osekowski}

\begin{equation} \label{Const Owsetsku}
Os(p) \stackrel{def}{=} 4 \sqrt{2} \cdot \left( \ \frac{p}{4} + 1 \ \right)^{1/p} \cdot \left( \  1 + \frac{p}{\ln(p/2)} \ \right), \ p \ge 4.
\end{equation}

 \ Notice that

$$
K  = K_{Os}  \stackrel{def}{=} \sup_{p \ge 4}  \left\{ \ \frac{Os(p)}{p/\ln p} \ \right\} \approx 15.7858...\ ,
$$

 so that

$$
Os(p) \le K_{Os} \cdot \frac{p}{\ln p}, \ p \ge 4.
$$

\vspace{3mm}

 \ It is proved by  Osekowski A, \ \cite{Osekowski},  that in the one - dimensional case

\vspace{3mm}

 \begin{equation} \label{Osek ineq}
 U(p;n) \le  Os(p) \cdot \max_{l = 1,2,\ldots,n} \mu(l,p), \ p \ge 4,
 \end{equation}

\vspace{3mm}

 \ following

\begin{equation} \label{K Osekov ineq}
 U(p;n) \le  K_{Os} \cdot \frac{p}{\ln p} \cdot \max_{l = 1,2,\ldots,n} \mu(l,p), \ p \ge 4.
 \end{equation}

\vspace{4mm}

 \hspace{3mm}  Our target in this preprint is extension the last estimate into the multidimensional, more precisely, matrix case. \par

 \vspace{5mm}

 \section{Moment estimations.}

 \vspace{4mm}

\hspace{3mm} Let the classical $ \ d \ - \ $  dimensional vector space $ \ X  = (X,||\cdot||) = R^d, \ d = 2,3,\ldots, \ $ be
 equipped with some  {\it arbitrary} finite for each vector $ \ x = \vec{x} \ $ norm $ \ ||x||, \ x = (x_1,x_2,\ldots,x_d) \in R^d. \ $
 Of course, $ \ ||x|| = 0 \ \Leftrightarrow x = 0. \ $ For instance,

\begin{equation} \label{lp norm}
||x||_p := \left[  \ \sum_{k=1}^d  |x_k|^p \ \right]^{1/p}, \ 1 \le p < \infty;
\end{equation}

  $ \ ||x||_{\infty} = \max_{k = 1,2,\ldots,d} |x_k|. \ $   As ordinary, the acting of the matrix $ \ A = \{ \ a(i,j) \ \}, \ i,j = 1,2,\ldots,d \ $
on the vector $ \  x = \vec{x} = \{x_j \} \ $  is defined as follows

$$
(Ax)_i = \sum_{j=1}^d a(i,j) \ x_j.
$$

 \vspace{3mm}

 \ The correspondent matrix norm for a square matrix $ \ A = \{ a_{i,j} \}, \ i,j = 1,2,\ldots,d \ $
will be denoted as follows. \par

\vspace{4mm}

 \ {\bf Definition 3.1.} \\

\vspace{3mm}

\begin{equation} \label{matrix norm}
||A|| \stackrel{def}{=} \sup_{0 \ne x \in R^d}  \left[ \ \frac{||Ax||}{||x||} \ \right].
\end{equation}

\vspace{3mm}

 \ For instance, let in (\ref{lp norm}) $ \ p = 2: \ $

\begin{equation} \label{l2 vec norm}
||x||_2 := \left[  \ \sum_{k=1}^d  |x_k|^2 \ \right]^{1/2}.
\end{equation}

 \ The correspondent matrix norm

\begin{equation} \label{spectral norm}
||A||_2 \stackrel{def}{=} \sup_{0 \ne x \in R^d}  \left[ \ \frac{||Ax||_2}{||x||_2} \ \right].
\end{equation}
is named as ordinary as {\it spectral norm} of the matrix $ \ A. \ $ It is known that

$$
||A||_2 = \sqrt{ \max \lambda(A^* A)},
$$
where $ \ \max \lambda(W) \ $ denotes a maximal (possible) eigenvalue of the symmetrical matrix $ \ W. \ $ \par

\vspace{4mm}

 \ Obviously, the set of all such a matrices (having a finite norm) forms a Banach algebra:

$$
||A \cdot B|| \le ||A|| \cdot ||B||.
$$

\vspace{3mm}

  \ The last definition  (\ref{matrix norm}) may be rewritten as follows. Denote by $ \ S = S(X) \ $ an unit sphere of this space:

$$
S  = S(X) \stackrel{def}{=} \{x, \ x \in R^d, \ ||x|| = 1 \}.
$$

 \ Then

\begin{equation} \label{spherical norm}
||A|| \stackrel{def}{=} \sup_{ x \in S} ||Ax||.
\end{equation}

 \ Moreover, denote by $ \  S_o = S_o(X)  \ $ the set of all {\it extreme points} of the unit  ball $ \ B(X) \ $
 of the space $ \ X = R^d: \  B(X) = \{x, x \in R^d, \ ||x|| \le 1 \ \}; \ $  then  $ \ S_o(X)  \subset S(X) \ $ and

\begin{equation} \label{extreme points}
||A|| \stackrel{def}{=} \sup_{ x \in S_o} ||Ax||  = \sup_{ x \in S_o(X)} ||Ax||,
\end{equation}
theorem of Krein - Milman   \cite{Krein Milman}.\par

 \ Furthermore, denote by $ \ Y = X^* \ $ the {\it conjugate space}  \ to source one  $ \ X = R^d. \ $  The acting the linear
 functional $ \ y \ $ from the space $ \ Y \ $ on the element $ \ x \ $ belonging to the first space $ \ X \ $  will be denoted
 as an inner product:

$$
(x,y) \stackrel{def}{=} y(x), \ x \in X, \ y \in Y.
$$

 \ The following relation holds, again by virtue of  Krein - Milman theorem

\begin{equation} \label{conj eq}
||x|| = \sup_{y \in S_o(Y)} (x,y) = \sup_{y \in S_o(X^*)} (x,y).
\end{equation}

 \vspace{5mm}

 \ Suppose now that the entries of the {\it sequence} of the
 square matrices   $ \ \Theta_n(i,j) \ $ of the size $ \ d \times d: \ i,j = 1,2,\ldots,d \ $  are {\it martingales}
 relative the index $ \ n \ $ and  the source filtration $ \ F(n) \ $ and with the correspondent martingale differences \ $  \zeta_k(i,j): \ $

 \begin{equation} \label{matr martdif}
 \Theta_n(i,j) = \sum_{k=1}^n \zeta_k(i,j).
 \end{equation}

\vspace{3mm}

 \hspace{3mm} Notice that

 \vspace{3mm}

\begin{equation} \label{norm reprsent}
 ||\Theta|| = \max_{\vec{x} \in S_0(X)} \ \max_{\vec{y} \in S_0(Y)} \ \left\{ \ \sum_{k = 1}^n \zeta_k(i,j) \ x_i \ y_j \ \right\},
\end{equation}

 \vspace{3mm}

or equally

 \vspace{3mm}

\begin{equation} \label{z norm reprsent}
 ||\Theta|| = \max_{{\vec{z} \in  S_0(X) \otimes S_0(Y)}} \left\{ \ \sum_{k = 1}^n \zeta_k(i,j) \ z(i,j)  \ \right\},
\end{equation}

 \vspace{3mm}

 where $ \ z = \vec{z} = \{ \ z(i,j) \ \}  = x_i \ y_j \in  Z \stackrel{def}{=}  S_0(X) \otimes S_0(Y). \ $ \par

 \vspace{3mm}

 \hspace{3mm} Let us consider the following {\it sequence } of the {\it naturally  normed}  random fields

\vspace{3mm}

 \begin{equation} \label{seq rf}
 L_n(z) \stackrel{def}{=} n^{-1/2}  \left\{ \ \sum_{k = 1}^n \zeta_k(i,j) \ z(i,j)  \ \right\}, \ z \in Z, \ n = 1,2,\ldots,
 \end{equation}

\vspace{3mm}

 \ so that $ \  L_n(z) = n^{-1/2}  \Theta(z) \ $ and following $ \ ||\Theta|| = \sqrt{n} \cdot ||L_n||.  \ $ \par

\vspace{3mm}

 \ We deduce  from the Osekovski estimate (\ref{K Osekov ineq})

\vspace{3mm}

\begin{equation} \label{Zz Osekov ineq}
| L_n(z)|_p  \le  K_{Os}\cdot \frac{p}{\ln p} \cdot \max_{l = 1,2,\ldots,n} \mu(l,p), \ p \ge 4,
\end{equation}

\vspace{3mm}

and quite analogously  $ \ | L_n(z_1)  - L_n(z_2)|_p  \le \ $

\vspace{3mm}

\begin{equation} \label{Z didderen  Osekov ineq}
 \le  K_{Os} \cdot \ ||z_1 - z_2|| \cdot \frac{p}{\ln p} \cdot \max_{l = 1,2,\ldots,n} \ \mu(l,p), \ p \ge 4, \ z_1,z_2 \in Z.
\end{equation}

 \ Denote by the sake of brevity

$$
\nu(p) :=  K_{Os} \cdot  \frac{p}{\ln p} \cdot \max_{l = 1,2,\ldots,n} \ \mu(l,p), \ p \ge 4;
$$

then

\vspace{3mm}

\begin{equation} \label{Zz  brev Osekov }
| L_n(z)|_p  \le \nu(p), \ p \ge 4, \ z \in Z.
\end{equation}

\vspace{3mm}

and

\vspace{3mm}

\begin{equation} \label{Zz dif brev}
 | L_n(z_1)  - L_n(z_2)|_p  \le ||z_1 - z_2|| \cdot \nu(p), \  \ p \ge 4, \ z_1,z_2 \in Z.
\end{equation}

\vspace{3mm}

 \ Denote by $ \ \kappa \ $  an entropic dimension of the set $ \ Z \ $ relative the classical Euclidean distance
  $ \  d(z_1, z_2) \stackrel{def}{=} |z_1 - z_2|. \ $   This means that

$$
H(Z, d, \epsilon) \le C + \kappa \ |\ln \epsilon|, \ \epsilon \in (0, 1/e).
$$

 \ Here  $ \ H(Z, d, \epsilon) \ $ denotes the {\it entropy} of the set $ \ Z \ $ at the point $ \ \epsilon, \ \epsilon > 0 \ $
 relative the distance function $ \ d; \ $ i.e. the natural logarithm  $ \ H(Z, d, \epsilon) := \ln N(Z, d, \epsilon) \ $
 of the number of the minimal $ \ d \ - \ \epsilon $ net of $ \ \epsilon - $
 balls in this distance $ \ d \ $ covering the whole set $ \ Z. \ $  Evidently,

\vspace{3mm}

\begin{equation} \label{dim estim}
 0 \le \kappa \le 2(d-1),
\end{equation}
and the case $ \ \kappa = 0 \ $ can not be excluded, still in the multidimensional case, when for instance $ \ B = [-1,1]^d. \ $ \par

\vspace{3mm}

 \hspace{3mm} One can apply the proposition of theorem 3.7.1 \ (4) of the monograph \cite{Ostrovsky1} to obtain the
 following estimate

\begin{equation} \label{key estim}
|\sup_{z \in Z} |L_n(z)| \ |_p \le \nu(p) \cdot \frac{10 p - \kappa}{p - \kappa}, \ p > \max(\kappa, 4).
\end{equation}

\vspace{3mm}

 \ Moreover, the random fields  $ \ L_n(z) \ $ are uniform  relative the index $ \ n \ $ uniform  continuous relative
 the parameter $ \ z; \ z \in Z. \ $ \par

\vspace{3mm}

  \hspace{3mm} As a consequence: define the variable

$$
\rho(p) :=  \nu(p) \cdot \frac{10 p - \kappa}{p - \kappa}, \ p > \max(\kappa, 4).
$$

\vspace{3mm}

 \hspace{3mm} {\sc  Let us impose henceforth the following important condition on the source random summands: }\par

\vspace{4mm}

\begin{equation} \label{key condit}
\exists  p_0 > \max(\kappa, 4) \ \Rightarrow \rho(p_0) < \infty.
\end{equation}

\vspace{3mm}

 \ Denote by $ \ p_- \ $ the {\it minimal} value for the finiteness of this function:

$$
p_- := \inf(p, \ p > \max(\kappa, 4), \  \rho(p) < \infty),
$$

as well as   by $ \ p_+ \ $ the {\it maximal} value for the finiteness of this function:

$$
p_+ := \sup(p, \ p > \max(\kappa, 4), \  \rho(p) < \infty).
$$
 \ The case when $ \ p_+ = \infty \ $ can not be excluded. \par

\vspace{3mm}

 \ Introduce correspondingly the following interval of finiteness of our function, i.e. such that

$$
J_0 := (p_-, \ p_+).
$$

\vspace{3mm}

 \hspace{3mm} {\bf Theorem 3.1.} We conclude under our notations and definitions

\vspace{3mm}

\begin{equation} \label{power estim}
| \ ||\Theta|| \ |_p \le \sqrt{n} \cdot  \ \rho(p), \ p \in J_0.
\end{equation}

\vspace{3mm}

\hspace{3mm} {\bf Remark 3.1.} Note that the last estimate (\ref{power estim}) is essentially, non - improvable still in the
one - dimensional case $ \ d = 1 \ $ and when the summands $ \ \{\xi(i)\} \ $ are centered and independent, see e.g. \cite{Ostrov6}. \par

 \vspace{3mm}

\hspace{3mm} {\bf Remark 3.2.}  The value  $ \ | \ ||\Theta|| \ |_q \ $ for the values $ \ q \in [1,  p_- ] $ may be estimates by
means of Lyapunov's inequality as follows. Let $  \   s \ $  be arbitrary point inside the interval  \ $ \ J_0: \ p_- < s < p_+; \ $
then

$$
| \ ||\Theta|| \ |_q  \le \sqrt{n} \cdot  \ \rho(s),
$$
therefore

$$
| \ ||\Theta|| \ |_q  \le \sqrt{n} \cdot  \  \inf_{s \in J_0} \rho(p).
$$

 \ Introduce  for simplicity the function

$$
\beta(p) := \rho(p), \ p \in J_0,
$$
 and

 $$
 \beta(p) := \inf_{s \in J_0} \rho(s), \  p \in [1,  p_-];
 $$
then the statement of theorem 3.1 may be rewritten  as follows

\vspace{3mm}

\begin{equation} \label{power complete estim}
| \ ||\Theta|| \ |_p \le \sqrt{n} \cdot  \ \beta(p), \ p \in [ 1, \ p_+)
\end{equation}

\vspace{4mm}

 \section{Exponential estimations. Main result.}

 \vspace{4mm}

 \hspace{3mm} We will deduce in this section  the {\it  exponential }  decreasing tail estimations  for the distribution the random
 variable $ \ ||\Theta||. \ $  We need to use for this purpose some facts from the the theory of the so - called Grand Lebesgue Spaces (GLS), which
 allow to deduce in particular the power and moreover the {\it exponential } decreasing estimates for tail of distribution for the random variables (r.v.),
 see e.g. \cite{Buldygin}, \cite{Ermakov}, \cite{Kozachenko1}, \cite{Liflyand}, \cite{Ostrovsky1}, \cite{FOS2022_Contemporary Mathematics}, \cite{ForKozOstr_Lithuanian}, etc. \par
 \hspace{3mm} Let $ \ (1,b), \ 1 < b \le \infty  \ $ be certain  numerical interval and let $ \ \psi = \psi(p), \ p \in (1,b) \ $ be some
 piece - continuous strictly positive numerical valued function. The so - called Grand Lebesgue Space (GLS) $ \ G \psi \ $ consists by definition
 on all the  numerical valued random variables $ \ \{\eta\} \ $ having a finite norm

\vspace{3mm}

\begin{equation} \label{norm Gpsi}
||\eta||G\psi \stackrel{def}{=} \sup_{ p \in (1,b)} \left\{ \ \frac{|\eta|_p}{\psi(p)} \ \right\}.
\end{equation}

\vspace{3mm}

 \ These spaces are (complete) rearrangement invariant (r.i.) functional Banach spaces. Further, define the new function
  $ \ h = h[\psi](p) = p \ln \psi(p) \ $ and its Young - Fenchel transform

$$
h^*(t) := \sup_{p \in (1,b)}(pt - h[\psi](p)), \ t \ge 1.
$$
 \ It is known that if $ \  0 < ||\eta||G\psi =: K < \infty,  \ $ then there holds a tail estimate, exponential in general case:

 \vspace{3mm}

\begin{equation} \label{tail est}
{\bf P}(|\eta| > t) \le 2 \exp \left(-h[\psi]^*( \ln(t/K)) \ \right), \ t \ge e K.
\end{equation}

\vspace{3mm}

 \ Moreover, the inverse conclusion also holds true, of course, under suitable natural conditions, \ \cite{Liflyand}, \cite{Kozachenko2}.\par

\vspace{4mm}

 \ {\bf Example 4.0.} \ Let $ \ m = \const > 0; \ $ the following propositions about the r.v. $ \ \eta \ $ are equivalent:

 $$
\sup_{p \ge 4} \left[ \ \frac{|\eta|_p}{p^{1/m}} \right] < \infty;
 $$

$$
\exists C = C(m) > 0 \ \Rightarrow \ {\bf P}  ( |\eta| > t) \le \exp \left(-C(m) t^m \right), \ t \ge 1.
$$

\ The case  $ \ m = 2 \ $  correspondent to the famous positive subgaussian random variable.

\vspace{3mm}

 \ {\bf Example 4.1.} \ Denote as usually by $ \  T[\xi](t)  \ $ the so - called {\it tail function} for the random variable $ \ \xi: \ $

$$
   T[\xi](t):= {\bf P}(|\xi| \ge t), \ t > 0. \
$$

\vspace{3mm}

   \ Let the r.v. $ \ \xi \ $ be such that
$$
T[\xi](t) \le  T^{(b, \gamma, S)}(t),  \ b = \const \in (1,\infty), \ \gamma = \const \in R,
$$
where by definition

$$
T^{(b, \gamma, S)}(t) \stackrel{def}{=} t^{-b} \ (\ln t)^{\gamma} S(\ln t), \ t \ge e,
$$
and

$$
\psi^{(b,\gamma,S)}(p) := (b - p)^{-(\gamma + 1)/b} \ S^{1/b} \left( \  \frac{1}{b - p}  \  \right), \ 1 \le p < b,
$$
where in turn   $ \ S = S(x), \ x \ge 1 \ $ is some positive continuous {\it slowly varying} function as $ \ x \to \infty. \ $ \par
  We have

$$
T_{\xi}(t) \le  T^{(b, \gamma, S)}(t) \ \Rightarrow ||\xi||G\psi^{(b,\gamma,S)} = C_1(b,\gamma,S) < \infty.
$$

 \ Inversely, if $ \ ||\xi||G\psi^{(b,\gamma,S)} = C_2 < \infty, \ $ then

$$
T[\xi](t) \le  C_3(b,\gamma,S) \cdot T^{(b, \gamma + 1, S)}(t);
$$
and both these estimates are non - improvable, see \cite{Kozachenko1} - \cite{Kozachenko2}. \par

\vspace{4mm}

 \hspace{3mm} As a consequence: \par

 \vspace{3mm}

 \hspace{3mm} {\bf Theorem 4.1.} We propose under  notations and conditions of theorem 3.1
 the following tail estimate, exponential in general case

\vspace{3mm}

\begin{equation} \label{main exp}
 {\bf P} ( n^{-1/2} \ ||\Theta|| > t)  \le  \exp \left( \ - h^*[\beta] (\ln t)  \ \right), \ t \ge e.
\end{equation}

\vspace{4mm}

 \hspace{3mm} {\bf Example  4.2.} Assume for instance that

 $$
 \max_{l =1,2,\ldots,n}\mu(l,p) \le C \ p^{\Delta} \ \ln p, \ p \in [4,\infty), \ C,\Delta = \const \in  (0,\infty);
 $$
 then

$$
 {\bf P} ( n^{-1/2} \ ||\Theta|| > t)  \le  C_2(\kappa,\Delta, C) \ t^{\kappa} \ \exp \left( \ - C_1(\kappa,\Delta) \ t^{1/(\Delta + 1)} \ \right), \ t \ge e,
$$

$ \ C_1 > 0, \ C_2 \in (0,\infty).\ $ \par

\vspace{4mm}

 \hspace{3mm} Let us illustrate the exactness of this estimate, by means of bringing of a (very simple) example.
We choose the one - dimensional "matrix" consisting from an unique centered random element $ \ \Theta = \xi, \ $ i.e. here $ \ n=d = 1. \ $
Impose the following condition  for these r.v.

\vspace{3mm}

$$
|\xi|_p \asymp p^{\Delta}, \ p \ge 1, \ \Delta = \const > 1.
$$

 \ There holds for this r.v.

$$
{\bf P}(|\xi| > t) \ge \exp \left(-C_0 \ t^{1/\Delta} \ \right), \ t \ge 1.
$$

\ One has in this example

$$
{\bf P} (|\Theta| > t) \ge  \exp \left(-C_0 \ t^{1/\Delta} \ \right), \ C_0 \in (0,\infty).
$$

\vspace{4mm}

 \hspace{3mm} {\bf Example  4.3.} Suppose now  that the r.v. $ \ \{\xi_{i,j} \}  \ $ are such that

$$
{\bf P} ( \ |\xi_{i,j}| > t \ ) \le t^{-b} \cdot (\ln t)^{\gamma} \cdot S(\ln t), \ t \ge e,
$$
where $ \ b = \const > 4, \ \gamma = \const \in R, \ S(\cdot) \ $ is as before  positive continuous slowly varying at infinity function.
There holds

$$
 {\bf P} ( n^{-1/2} \ ||\Theta|| > t)  \le C_4(b, \gamma, S(\cdot), d) \cdot t^{-b} \cdot (\ln t)^{\gamma + 1} \cdot S(\ln t), \ t \ge e.
$$

\vspace{4mm}

 \hspace{3mm} Let us illustrate again the exactness of this estimate, by means of bringing of a (very simple) example.
We choose as before the one - dimensional "matrix" consisting from an unique centered random element $ \ \Theta = \xi, \ $ i.e. here $ \ n=d = 1. \ $
Impose the following distribution for these r.v.

$$
{\bf P}(|\xi| > t) := t^{-b} \ (\ln t)^{\gamma} \ S(\ln t)  = T^{b,\gamma,S}(t), \ t \ge t_0 = \const > e.
$$

\ One has in this example

$$
{\bf P} (|\Theta| > t) \ge  T^{b,\gamma,S}(t) = t^{-b} \ (\ln t)^{\gamma} \ S(\ln t), \ t \ge t_0 .
$$

\vspace{4mm}

\section{Concluding remarks.}

\vspace{3mm}

 \hspace{3mm} It is no hard  to generalize the obtained results into the non - square martingale matrix.
  A more interesting problem:  to extend  these estimates into infinite - dimensional martingale linear operators.\par

 \vspace{6mm}

\emph{\textbf{\footnotesize Acknowledgements}.} {\footnotesize M.R. Formica is member of Gruppo
Nazionale per l'Analisi Matematica, la Probabilit\`{a} e le loro Applicazioni (GNAMPA) of the Istituto Nazionale di Alta Matematica (INdAM) and member of the UMI group \lq\lq Teoria dell'Approssimazione e Applicazioni (T.A.A.)\rq\rq and is partially supported by the INdAM-GNAMPA project, {\it Risultati di regolarit\`{a} per PDEs in spazi di funzione non-standard}, codice CUP\_E53C22001930001.
}

\end{document}